\documentclass[12pt]{article}
\usepackage{amsfonts}
\usepackage{amssymb}
\usepackage{amsbsy}
\usepackage{amsthm}

\newtheorem{thm}{Theorem}[section]
\newtheorem{cor}[thm]{Corollary}
\newtheorem{lem}[thm]{Lemma}
\newtheorem{prop}[thm]{Proposition}
\theoremstyle{definition}
\newtheorem{defn}[thm]{Definition}
\theoremstyle{remark}
\newtheorem{rem}[thm]{Remark}

\newcommand{\bt}{\begin{thm}}
\newcommand{\et}{\end{thm}}
\newcommand{\bc}{\begin{cor}}
\newcommand{\ec}{\end{cor}}
\newcommand{\bl}{\begin{lem}}
\newcommand{\el}{\end{lem}}
\newcommand{\bp}{\begin{prop}}
\newcommand{\ep}{\end{prop}}
\newcommand{\bd}{\begin{defn}}
\newcommand{\ed}{\end{defn}}
\newcommand{\br}{\begin{rem}}
\newcommand{\er}{\end{rem}}
\newcommand{\bpr}{\begin{proof}}
\newcommand{\epr}{\end{proof}}

\newcommand{\bi}{\begin{itemize}}
\newcommand{\ei}{\end{itemize}}
\newcommand{\be}{\begin{enumerate}}
\newcommand{\ee}{\end{enumerate}}
\newcommand{\ds}{\displaystyle}
\newcommand{\ba}{\begin{array}}
\newcommand{\ea}{\end{array}}
\newcommand{\beq}{\begin{equation}}
\newcommand{\eeq}{\end{equation}}
\newcommand{\beqa}{\begin{eqnarray}}
\newcommand{\eeqa}{\end{eqnarray}}

\newcommand{\N}{{\mathbb N}}
\newcommand{\Z}{{\mathbb Z}}
\newcommand{\R}{{\mathbb R}}
\newcommand{\C}{{\mathbb C}}
\newcommand{\T}{{\mathbb T}}
\newcommand{\D}{{\mathbb D}}
\newcommand{\PP}{{\mathbb P}}

\newcommand{\cA}{{\mathcal  A}}

\newcommand{\cC}{{\mathcal  C}}

\newcommand{\cH}{{\mathcal  H}}

\newcommand{\cX}{{\mathcal  X}}
\newcommand{\cS}{{\mathcal  S}}

\newcommand{\bfa}{{\mathbf a}}
\newcommand{\bfb}{{\mathbf b}}

\newcommand{\bfX}{{\mathbf X}}

\newcommand{\bsF}{{\boldsymbol \Phi}}

\newcommand{\spn}{\mathrm{span}}
\newcommand{\supp}{\mathrm{supp}}
\newcommand{\spec}{\sigma}

\newcommand{\limsn}{\mathop{\overline{\lim}} \limits_n}

\newcommand{\re}{\mathrm{Re}}
\newcommand{\im}{\mathrm{Im}}

\begin{document}

\title{\bf Minimal representations of unitary operators and orthogonal polynomials
on the unit circle
\footnote{This work was supported by Project E-12/25 of DGA (Diputaci\'on General de
Arag\'on) and by Ibercaja under grant IBE2002-CIEN-07.}}
\author{M.J. Cantero, L. Moral
\footnote{Corresponding author. Tel.: +34-976-76-11-41. Fax: +34-976-76-11-25.
E-mail address: lmoral@unizar.es}
, L. Vel\'azquez}
\date{\small \centerline{Departamento de Matem\'atica Aplicada, Universidad de Zaragoza,}
\hfill\break\vskip-25pt
\centerline{Pza. San Francisco s/n, 50009 Zaragoza, Spain}}
\maketitle

\begin{abstract}
In this paper we prove that the simplest band representations of unitary operators on a
Hilbert space are five-diagonal. Orthogonal polynomials on the unit circle play an
essential role in the development of this result, and also provide a parametrization of
such five-diagonal representations which shows specially simple and interesting
decomposition and factorization properties. As an application we get the reduction of
the spectral problem of any unitary Hessenberg matrix to the spectral problem of
a five-diagonal one. Two applications of these results to the study of orthogonal
polynomials on the unit circle are presented: the first one concerns Krein's Theorem;
the second one deals with the movement of mass points of the orthogonality measure under
monoparametric perturbations of the Schur parameters.
\end{abstract}


\medskip

\noindent{\it (2000) AMS Mathematics Subject Classification}: 42C05, 47B36.

\section{Introduction}

Matrix representations are an important tool for the study of linear operators on a
Hilbert space. They allow, for instance, the use of perturbation techniques for the
comparison of operators defined on different Hilbert spaces. Besides, the freedom in the
choice of the representation can be used to get a simple one that can make the analysis
of the operator easier. Usually a band representation with minimum size band is
desirable. A band matrix $(c_{i,j})$ is $(p,q)$-diagonal if $c_{i,j}=0$ for $i-j< p$ and
$j-i > q$. A matrix that is $([{n \over 2}],[{n-1 \over 2}])$-diagonal or $([{n-1 \over
2}],[{n \over 2}])$-diagonal, is called a $n$-diagonal matrix. If every operator of a
certain class has a $n$-diagonal representation but not all of them have a $n-1$-diagonal
one, we say that $n$-diagonal representations are the minimal representations of the
class.

Concerning the class of self-adjoint operators, any two-diagonal representation must be
diagonal due to its symmetry, but a diagonal representation is only possible in the case
of pure point spectrum. Therefore, the minimal representations are at least tri-diagonal.
In fact, they are tri-diagonal since, as a consequence of the spectral theorem, every
self-adjoint operator is unitarily equivalent to an orthogonal sum of self-adjoint
multiplication operators \cite{ReSi72} and, hence, the use of basis of orthogonal
polynomials on the real line gives a tri-diagonal representation \cite{St32}.

Unitary operators, together with self-adjoint ones, are the most important examples of
normal operators. However, in spite of their importance, the minimal representations for
unitary operators are an open problem. Analogously to the self-adjoint case, the study
can be reduced to unitary multiplication operators, but the use of basis of orthogonal
polynomials on the unit circle then leads to Hessenberg instead of band representations
\cite{Ge44,Al77,Gr93,MaGo91,Te92}. As for the possibility of band representations, it has
been recently proved in \cite{BoHoJo} that any unitary tri-diagonal matrix decomposes as
a sum of $1\times1$ and $2\times2$ diagonal blocks and, therefore, it has a pure point
spectrum. This shows that the minimal representations of unitary operators are at least
four-diagonal. B. Simon has conjectured in a preliminary version of \cite{Si04} that a
similar decomposition should happen for any unitary four-diagonal matrix, which would
imply that the minimal representations for unitary operators are at least five-diagonal.

A step to get the minimal representations of unitary operators was taken by the
authors in \cite{FIVE}. The results presented there imply that any unitary operator has a
five-diagonal representation. In the next section we introduce these five-diagonal
representations and their connection with orthogonal polynomials on the unit circle.
Section 3 is devoted to the study of such representations and their properties.
Although Hessenberg matrices have been extensively studied, all this analysis will be
done jointly for five-diagonal and Hessenberg representations for several reasons:

\noindent -- It is convenient to understand the improvements given by the five-diagonal
representations, if compared with the known Hessenberg ones. Some concrete examples of
the advantages of the five-diagonal representations will be clearly shown in the
applications discussed in sections 4 and 5.

\noindent -- The connections between Hessenberg and five-diagonal representations provide
an algorithm that reduces the spectral problem of any unitary Hessenberg matrix to the
spectral problem of a five-diagonal one (``five-diagonal reduction" of the spectral
problem of a unitary Hessenberg matrix). The importance of this result is due to the
increasing interest in the study of unitary Hessenberg matrices in numerical linear
algebra \cite{Gr93,Gr86,GrRe90} and digital signal processing applications (see
\cite{DeGe90} and references therein).

\noindent -- The analysis of unitary Hessenberg matrices is the main tool to prove that
the minimal representations of unitary operators are indeed five-diagonal. This result is
a consequence of a more general one that closes Section 3: $(1,q)$-diagonal or
$(p,1)$-diagonal representations of unitary operators are possible only in the case of
pure point spectrum.

In sections 4 and 5 we consider some applications of the minimal representations of
unitary operators to the study of orthogonal polynomials on the unit circle. Both
applications concern the relation between the support of the measure of orthogonality and
the corresponding Schur parameters. Section 4 shows the advantages of the five-diagonal
representation for the analysis of the limit points of the support of the measure, while
Section 5 is devoted to the study of the isolated mass points. We finish this last
section giving several explicit examples of perturbations of the Schur parameters that
keep an arbitrary mass point invariant.

Now we proceed with the conventions for the notation. For any subset $\cA$ of a Hilbert
space, $\overline \cA$ is its closure and $\spn \cA$ the set of all finite linear
combinations of $\cA$. Also, if $\cS$ is a subspace of the Hilbert space, $\cA^{\bot\cS}$
means the subspace of $\cS$ orthogonal to $\cA$.

Given a linear operator $T$ on a Hilbert space, $T^*$ denotes its adjoint and $\spec(T)$
its spectrum, while for every complex matrix $M$, $M^T$ is its transpose and
$M^* = \overline M^T$. $I$ and $I_N$ represent the unit matrix of order infinite and $N$,
respectively. Any matrix of order $N$ is considered as an operator in $\C^N$, and any
infinite bounded matrix is identified with the continuous operator that it defines in
$\ell^2$, the Hilbert space of square-sumable sequences in $\C$.  The inner products in
$\C^N$ and $\ell^2$ are denoted by $(\cdot,\cdot)$, and the corresponding canonical basis
by $\{e_n\}$. No misunderstanding will arise from this common notation.

The term measure always means non-negative finite Borel measure, and, without loss of
generality, we will consider only probability measures. If $\mu$ is a measure on a
subset of $\C$, $\supp\mu$ is its support and $L^2_\mu$ the Hilbert space of
$\mu$-square-integrable complex functions with inner product
$$
\big< f,g \big>_\mu := \int f(z) \overline{g(z)} \, d\mu(z),
\quad f,g \in L^2_\mu.
$$
$\T:=\{z\in\C:|z|=1\}$ is the unit circle and $\D:=\{z\in\C:|z|<1\}$ the open unit disk
in the complex plane. A multiplication operator on $\T$ has the form
$$
U_\mu \colon \mathop{L^2_\mu \to L^2_\mu} \limits_{f(z) \; \to \; zf(z)}
$$
where $\mu$ is a measure on $\T$.

\section{Representations of unitary operators and orthogonal polynomials on $\T$}

Given a unitary operator U on a separable Hilbert space $\cH$, the equivalence between
the following assertions is known \cite{St32}:
\bi
\item The spectrum of $U$ is simple.
\item $U$ has a cyclic vector $v \in \cH$, in the sense of
      $\overline{\spn\{U^nv\}}_{n\in\Z}=\cH$.
\item $U$ is unitarily equivalent to a multiplication operator on $\T$.
\ei

A standard application of Zorn's lemma shows that any unitary operator can be expressed
as a (finite or infinite) orthogonal sum of unitary operators with cyclic vectors.
Therefore, the study of unitary operators becomes the study of multiplication operators
on $\T$. As for the spectral properties of such multiplication operators, it is known
that $\spec(U_\mu)=\supp\mu$, the mass points of $\mu$ being the eigenvalues of $U_\mu$.
The eigenvectors of $U_\mu$ associated with an eigenvalue $\lambda$ are spanned by the
characteristic function $\cX_\lambda$ of the set $\{\lambda\}$.

For a long time, the usual attempts to get matrix representations of $U_\mu$ have dealt
with basis constituted by orthogonal polynomials (OP) with respect to $\mu$, that is,
polynomials satisfying \beq \label{OP} \deg\varphi_n=n, \quad \big< \varphi_n, \varphi_m
\big>_\mu = \delta_{n,m}, \quad n,m\geq0. \eeq When $\supp\mu$ has a finite number $N$ of
elements, dim$(L^2_\mu)=N$ and such a basis $\bsF_N:=(\varphi_n)_{n=0}^{N-1}$ comes from
the orthogonalization of $\{z^n\}_{n=0}^{N-1}$. $\bsF_N$ is called a finite segment of OP
associated with $\mu$. If $\supp\mu$ is infinite, dim$(L^2_\mu)=\aleph_0$ and the
orthogonalization of the infinite set $\{z^n\}_{n\ge0}$ gives a sequence
$\bsF:=(\varphi_n)_{n\ge0}$ satisfying (\ref{OP}) that is called a sequence of OP with
respect to $\mu$. However, such a sequence is not always a basis of $L^2_\mu$ since the
polynomials are not always dense in $L^2_\mu$.

In what follows, $\varphi_n$ denotes the unique $n$-th OP with respect to $\mu$ with
positive leading coefficient $\kappa_n$. It is known that these polynomials satisfy the
recurrence relation
\beq \label{RR-OP}
\ba{l}
\varphi_0(z) = 1,
\\
\rho_n \varphi_n(z) = z\varphi_{n-1}(z) + a_n \varphi_{n-1}^*(z), \quad n\geq1,
\ea
\eeq
where $p^*(z) := z^n \overline p(z^{-1})$ for a polynomial $p$ of degree $n$,
$\rho_n:=\sqrt{1-|a_n|^2}$ and $a_n\in\D$ are known as the Schur parameters associated
with $\mu$.

Besides, when $\supp\mu=\{z_1,z_2,\dots,z_N\}$, the same arguments that give
(\ref{RR-OP}) show that the polynomial $\psi(z) = (z-z_1)(z-z_2)\cdots(z-z_N)$ satisfies
\beq \label{POP}
\kappa_{N-1} \psi(z) = z\varphi_{N-1}(z) + a_N \varphi_{N-1}^*(z), \quad a_N\in\T.
\eeq

It is known that the preceding results establish a one to one correspondence between:
\bi
\item Probability measures supported on $N$ points of the unit circle and vectors
      $\bfa_N:=(a_1,a_2,\dots,a_N) \in \D^{N-1}\times\T$.
\item Probability measures supported on an infinite subset of the unit circle and
      sequences $\bfa:=(a_n)_{n\in\N} \in \D^{\aleph_0}$.
\ei

If $\bsF$ is the sequence of OP related to a measure $\mu$ with infinite support, from
(\ref{RR-OP}) we find that the matrix of $U_\mu$ with respect to $\bsF$ is a Hessenberg
one given by \cite{Ge44,Al77,Gr93,MaGo91,Te92}
$$
H(\bfa):=\pmatrix{
-a_1 & \kern-3pt -\rho_1 a_2 & \kern-3pt -\rho_1 \rho_2 a_3 &
\kern-3pt -\rho_1 \rho_2 \rho_3 a_4 & \kern-3pt \cdots
\cr
\kern3pt \rho_1 & \kern-3pt -\overline a_1 a_2 &
\kern-3pt -\overline a_1 \rho_2 a_3 &
\kern-3pt -\overline a_1 \rho_2 \rho_3 a_4 & \kern-3pt \cdots
\cr
 & \rho_2 & \kern-3pt -\overline a_2 a_3 &
\kern-3pt -\overline a_2 \rho_3 a_4 & \kern-3pt \cdots
\cr
 & & \rho_3 & \kern-3pt -\overline a_3 a_4 & \kern-3pt \cdots &
\cr
 & & & \rho_4 & \kern-3pt \cdots
\cr
 & & & & \kern-3pt \cdots\cr
},
$$
where $\bfa$ is the corresponding sequence of Schur parameters.

The principal matrix of order $N$ of $H(\bfa)$ only depends on the vector $\bfa_N$ and
will be denoted by $H(\bfa_N)$. If $\bfa_N\in\D^{N-1}\times\T$ and $\mu$ is the
related finitely supported measure we get from (\ref{RR-OP}) and (\ref{POP}) that
$z\bsF_N(z) = H(\bfa_N)^T\bsF_N(z) + \kappa_{N-1}\psi(z)e_N$, $\bsF_N$ being
the corresponding finite segment of OP \cite{MPOP}. Since $\psi(z)=0$ $\mu$-a.e.,
we find that $H(\bfa_N)$ is the matrix of $U_\mu$ with respect to $\bsF_N$.

Apart from its complexity, the infinite matrix $H(\bfa)$ represents the full operator
$U_\mu$ only when the polynomials are dense in $L^2_\mu$, that is, when
$\bfa\not\in\ell^2$ \cite{Ge61,Sz75}. In general, $H(\bfa)$ represents the restriction of
$U_\mu$ to the closure of $\PP:=\spn\{z^n\}_{n\ge0}$. Hence, although $H(\bfa)$ is
always isometric, it is unitary iff $z\overline \PP = \overline \PP$. Since this
condition is equivalent to $\overline\PP=L^2_\mu$, we see that $H(\bfa)$ is unitary iff
$\bfa\not\in\ell^2$.

The measures corresponding to sequences $\bfa\in\ell^2$ constitute the so-called
Szeg\H o class. A possibility of getting a matrix representation for $U_\mu$ in this case
is to enlarge the OP basis to get an orthonormal basis of $L^2_\mu$. This possibility
is exploited in \cite{Si04}, obtaining a doubly infinite unitary matrix in which
$H(\bfa)$ is embedded. Anyway, the complexity of the matrix representation remains.

If we want to simplify the matrix representation of $U_\mu$ solving at the same time the
problem for the Szeg\H o class, we have to change completely the choice of the basis for
$L^2_\mu$. Since the space of Laurent polynomials is always dense in $L^2_\mu$, a more
natural choice for a basis is the orthogonal Laurent polynomials (OLP) with respect to
$\mu$, related to the corresponding OP by \cite{FIVE,Th88}
\beq \label{OP-OLP}
\chi_{2k}(z) = z^{-k} \varphi_{2k}^*(z), \quad \chi_{2k+1} = z^{-k} \varphi_{2k+1}(z),
\quad k\geq0.
\eeq
The above relation gives a finite segment of OLP $\bfX_N:=(\chi_n)_{n=0}^{N-1}$ in the
case of a measure supported on $N$ points, or a sequence $\bfX:=(\chi_n)_{n\ge0}$ of OLP
for an infinitely supported measure. $\bfX_N$ and $\bfX$ always constitute an
orthonormal basis of the corresponding space $L^2_\mu$.

If the measure $\mu$ has infinite support, we get from (\ref{RR-OP}) the following matrix
representation for the operator $U_\mu$ with respect to the related sequence $\bfX$ of
OLP \cite{FIVE}
$$
C(\bfa):=\pmatrix{
-a_1 & \kern-7pt -\rho_1 a_2 & \rho_1 \rho_2
\cr
\kern7pt \rho_1 & \kern-7pt -\overline a_1 a_2 & \overline a_1 \rho_2 & 0
\cr
\kern5pt 0 & \kern-7pt -\rho_2 a_3 & \kern-7pt -\overline a_2 a_3 &
\kern-7pt -\rho_3 a_4 & \rho_3 \rho_4
\cr
& \rho_2 \rho_3 & \overline a_2 \rho_3 & \kern-7pt -\overline a_3 a_4  &
\overline a_3 \rho_4 & 0
\cr
&& 0 & \kern-7pt -\rho_4 a_5 & \kern-7pt -\overline a_4 a_5 &
\kern-7pt -\rho_5 a_6 & \rho_5 \rho_6
\cr
&&& \rho_4 \rho_5 & \overline a_4 \rho_5 & \kern-7pt -\overline a_5 a_6 &
\overline a_5 \rho_6 & 0
\cr
&&&& \hskip-35pt\ddots & \hskip-35pt\ddots &
\hskip-35pt\ddots & \hskip-20pt\ddots & \ddots},
$$
$\bfa$ being the corresponding sequence of Schur parameters.

Now we deal with a five-diagonal matrix that, contrary to the Hessenberg one,
always represents the full operator $U_\mu$ and, hence, is unitary for any
$\bfa\in\D^{\aleph_0}$. Besides, it has a much simpler dependence on the Schur
parameters.

The principal matrix of order $N$ of $C(\bfa)$, that only depends on $\bfa_N$, will be
denoted by $C(\bfa_N)$. Analogously to the case of the Hessenberg representation, if
$\bfa_N\in\D^{N-1}\times\T$ and $\mu$ is the related measure, $C(\bfa_N)$ is the matrix
of $U_\mu$ with respect to the corresponding finite segment of OLP $\bfX_N$: using
(\ref{RR-OP}), (\ref{POP}) and (\ref{OP-OLP}) we find that
$z\bfX_N(z) = C(\bfa_N)^T \bfX_N(z) + \bfb_N z^{-[{N-1\over2}]} \psi(z)$, where
$$
\bfb_N=\cases{
\kappa_{N-1} e_N, & if $N$ is even,
\cr
\kappa_{N-1} (\rho_{N-1}e_{N-1}+\overline a_{N-1}e_N), & if $N$ is odd,
}
$$
and, since $\psi(z)=0$ $\mu$-a.e., we get the desired result.

\medskip

Let $\mu$ be a measure on $\T$, $(\varphi_n)$ the corresponding OP and $(\chi_n)$ the
related OLP. As a consequence of the whole previous discussion, if $\mu$ is associated
with the sequence $\bfa\in\D^{\aleph_0}$ of Schur parameters, $\sigma(H(\bfa))=\supp\mu$
for $\bfa\notin\ell^2$ while $\sigma(C(\bfa))=\supp\mu$ always happens. Also, if $\mu$
is the finitely supported measure associated with $\bfa_N\in\D^{N-1}\times\T$, then
$\sigma(H(\bfa_N))=\sigma(C(\bfa_N))=\supp\mu$. Similar relations hold between the mass
points of the measure and the eigenvalues of the related matrices.
As for the eigenvectors associated with a mass point $\lambda$, since
$\big<\cX_\lambda,\varphi_n\big>_\mu=\mu(\{\lambda\})\overline{\varphi_n(\lambda)}$
and $\big<\cX_\lambda,\chi_n\big>_\mu=\mu(\{\lambda\})\overline{\chi_n(\lambda)}$, we
find that $\sum_n\overline{\varphi_n(\lambda)}e_n$ is an eigenvector of the corresponding
Hessenberg matrix when it represents the full operator $U_\mu$, while
$\sum_n\overline{\chi_n(\lambda)}e_n$ is always an eigenvector of the related
five-diagonal matrix.

Let $\lambda$ be a mass point of $\mu$. Using the decomposition of $\cX_\lambda$ with
respect to the OLP basis we find that $\mu(\{\lambda\}) =
\big<\cX_\lambda,\cX_\lambda\big>_\mu = \mu(\{\lambda\})^2 \sum_n |\chi_n(\lambda)|^2$
and, so, since $\lambda\in\T$, we get from the relation (\ref{OP-OLP}) between OP and OLP
that $\mu(\{\lambda\}) = 1/\sum_n |\varphi_n(\lambda)|^2$. Thus, when $\mu$ has infinite
support, $\bsF(\lambda)\in\ell^2$ if $\lambda$ is a mass point. Conversely, if
$\lambda\in\T$ is such that $\bsF(\lambda)\in\ell^2$, then $\lambda$ is a mass point
since $C(\bfa)^*\overline{\bfX(z)} = \overline z \overline{\bfX(z)}$, $\forall z\in\C$.
Notice that these arguments also work using the OP basis but restricted to measures
outside the Szeg\H o class.

\medskip

Among other things, the preceding results show that the minimal representations of
unitary operators are at most five-diagonal, but, are they exactly five-diagonal?

Moreover, like any unitary operator, every Hessenberg matrix that is unitary must be
unitarily equivalent to a five-diagonal one. However, a question remains if we want
to complete the ``five-diagonal reduction'' of the spectral problem for any unitary
Hessenberg matrix: which one is the five-diagonal matrix related to an arbitrary unitary
Hessenberg one?

A deeper study of unitary five-diagonal and Hessenberg matrices will answer the above
questions.

\section{Five-diagonal and Hessenberg matrices}

The five-diagonal matrices presented in the previous section are examples of the
following kind of matrices, that can be considered an intermediate step between the
five-diagonal and the tri-diagonal case.

\bd \label{PARA-TRI}
A (finite or infinite) five-diagonal matrix $C=\left(c_{i,j}\right)$ is called
para-tridiagonal if $c_{2k,2k+2}=c_{2k+1,2k-1}=0$, $\forall k\ge1$, that is,
$$
C = \pmatrix{
c_{1,1} & c_{1,2} & c_{1,3} \cr
c_{2,1} & c_{2,2} & c_{2,3} & 0  \cr
0 & c_{3,2} & c_{3,3} & c_{3,4} & c_{3,5} \cr
& c_{4,2} & c_{4,3} & c_{4,4} & c_{4,5} & 0 \cr
& & 0 & c_{5,4} & c_{5,5} & c_{5,6} & c_{5,7} \cr
& & & c_{6,4} & c_{6,5} & c_{6,6} & c_{6,7} & 0 \cr
& & & & \ddots & \ddots & \ddots & \ddots & \ddots
}.
$$
$C$ is called irreducible if $c_{2,1}\neq0$ and $c_{2k-1,2k+1},c_{2k+2,2k}\neq0$,
$\forall k\geq1$.
\ed

Unitary irreducible para-tridiagonal matrices, like unitary irreducible Hessenberg ones,
have $e_1$ as a cyclic vector (in fact, any vector $e_n$, $n\in\N$, is cyclic in the
first case). Therefore, any unitary irreducible para-tridiagonal matrix is the matrix
representation of a multiplication operator on $\T$ and, hence, is unitarily equivalent
to one with the form $C(\bfa)$ or $C(\bfa_N)$. However we do not know how to describe
this relation exactly. The following result is the first step to answer this and the
previous questions, since it provides the general form of infinite unitary
para-tridiagonal and Hessenberg matrices. The matrix representations introduced in the
preceding section are indispensable guides for taking this step.

\bt \label{INFINITE}
An infinite para-tridiagonal (Hessenberg) matrix is unitary (isometric) iff it has the
form $C(\bfa,\bfb)$ ($H(\bfa,\bfb)$), where
$\bfa,\bfb \in \C^{\aleph_0}$ are such that $|a_n|^2+|b_n|^2=1$, $\forall n\in\N$, and
$$
\ba{l}
C(\bfa,\bfb) := \pmatrix{
-a_1 & \kern-7pt -\overline b_1 a_2 &
\overline b_1 b_2
\cr
\kern7pt b_1 & \kern-7pt -\overline a_1 a_2 &
\overline a_1 b_2 & 0
\cr
\kern5pt 0 & \kern-7pt -\overline b_2 a_3 & \kern-7pt -\overline a_2 a_3 &
\kern-7pt -\overline b_3 a_4 & \overline b_3 b_4
\cr
& \overline b_2 b_3 & \overline a_2 b_3 & \kern-7pt -\overline a_3 a_4  &
\overline a_3 b_4 & 0
\cr
&& 0 & \kern-7pt -\overline b_4 a_5 & \kern-7pt -\overline a_4 a_5 &
\kern-7pt -\overline b_5 a_6 & \overline b_5 b_6
\cr
&&& \overline b_4 b_5 & \overline a_4 b_5 & \kern-7pt -\overline a_5 a_6 &
\overline a_5 b_6 & 0
\cr
&&&& \hskip-35pt\ddots & \hskip-35pt\ddots &
\hskip-35pt\ddots & \hskip-20pt\ddots & \ddots},
\medskip\\
H(\bfa,\bfb):=\pmatrix{
-a_1 & \kern-3pt -\overline b_1 a_2 & \kern-3pt -\overline b_1 \overline b_2 a_3 &
\kern-3pt -\overline b_1 \overline b_2 \overline b_3 a_4 & \kern-3pt \cdots
\cr
\kern3pt b_1 & \kern-3pt -\overline a_1 a_2 &
\kern-3pt -\overline a_1 \overline b_2 a_3 &
\kern-3pt -\overline a_1 \overline b_2 \overline b_3 a_4 & \kern-3pt \cdots
\cr
 & b_2 & \kern-3pt -\overline a_2 a_3 &
\kern-3pt -\overline a_2 \overline b_3 a_4 & \kern-3pt \cdots
\cr
 & & b_3 & \kern-3pt -\overline a_3 a_4 & \kern-3pt \cdots &
\cr
 & & & b_4 & \kern-3pt \cdots
\cr
 & & & & \kern-3pt \cdots\cr
}.
\ea
$$
\et

\bpr
An infinite para-tridiagonal matrix $C$ can be written in the way
$$
C=\pmatrix{
C_1^T & C_2 & 0 & \cdots \cr
0 & C_3^T & C_4 & \cdots \cr
0 & 0 & C_5^T & \cdots \cr
\cdots & \cdots & \cdots & \cdots
},
\quad \ba{l} C_1\in\C^{(1,2)}, \\ C_n\in\C^{(2,2)}, \; n\geq2. \ea
$$
It is unitary iff $C^*C=CC^*=I$, which is equivalent to
$$
\ba{l}
C_1C_1^*=I_1, \\
C_nC_n^*+(C_{n-1}^*C_{n-1})^T=I_2, \quad n\geq2, \\
C_n^*C_{n-1}^T=0, \quad n\geq2.
\ea
$$
The first condition means that
$$
C_1=\pmatrix{-a_1 & b_1}, \quad |a_1|^2+|b_1|^2=1,
$$
and, then, by induction, we find that the rest of the equations are satisfied iff
$$
C_n = \pmatrix{
-\overline b_{n-1} a_n & \overline b_{n-1} b_n \cr
-\overline a_{n-1} a_n & \overline a_{n-1} b_n},
\quad |a_n|^2+|b_n|^2=1, \quad n\geq2.
$$
This proves the theorem in the para-tridiagonal case.

Now, let $H$ be a Hessenberg matrix, that is, its $n$-th column $h_n$ belongs to
$\spn\{e_1,e_2,\dots,e_{n+1}\}\subset\ell^2$. $H$ is isometric iff $H^*H=I$, which means
that $\{h_n\}_{n\in\N}$ is an orthonormal set of $\ell^2$. We will see that this is
equivalent to
\beq \label{COL}
\ba{l}
h_n = b_n e_{n+1} - a_n v_n, \quad n\in\N
\\
v_n = \sum_{i=1}^n
\overline a_{i-1} \overline b_i \overline b_{i+1} \cdots \overline b_{n-1} e_i,
\quad n\in\N,
\\
a_0=1; \quad |a_n|^2+|b_n|^2=1, \quad n\in\N,
\ea
\eeq
which proves the theorem for the Hessenberg case.

First of all, let us suppose that the columns of $H$ have the form (\ref{COL}).
From the expression of $v_n$ we find that
$v_{n+1} = \overline b_n v_n + \overline a_n e_{n+1}$
for $n\in\N$. Therefore, we get by induction that
$v_n \bot \{h_1,h_2,\dots,h_{n-1}\}$ and $(v_n,v_n)=1$ for $n\in\N$.
Then, the expression for $h_n$ implies that
$h_n \bot \{h_1,h_2,\dots,h_{n-1}\}$ and $(h_n,h_n)=1$ for $n\in\N$.

On the other hand, if the columns of $H$ form an orthonormal set of $\ell^2$, then
we can write $h_1=-a_1e_1+b_1e_2$, $|a_1|^2+|b_1|^2=1$, and, for $n\geq2$,
$h_n=b_ne_{n+1}+u_n$, $b_n\in\C$,
$u_n\in\{h_1,h_2,\dots,h_{n-1}\}^{\bot\spn\{e_1,e_2,\dots,e_n\}}$.
So, $h_1$ has the form given by (\ref{COL}). Moreover, let us suppose that
$h_1,h_2,\dots,h_{n-1}$ satisfy (\ref{COL}). Then,
$\{h_1,h_2,\dots,h_{n-1}\}^{\bot\spn\{e_1,e_2,\dots,e_n\}}=\spn\{v_n\}$
and we find that $u_n=-a_nv_n$, $a_n\in\C$. The condition $(h_n,h_n)=1$ gives
$|a_n|^2+|b_n|^2=1$. This proves by induction that $h_n$ has the form ($\ref{COL}$)
for $n\in\N$.
\epr

A consequence of Theorem \ref{INFINITE} is its analogue for finite matrices. The result
for the Hessenberg case was already known \cite{Gr86}. In what follows, since the
principal submatrix of order $N$ of $C(\bfa,\bfb)$ ($H(\bfa,\bfb)$) only depends on
$\bfa_N,\bfb_{N-1}$, this submatrix will be denoted by $C(\bfa_N,\bfb_{N-1})$
($H(\bfa_N,\bfb_{N-1})$).

\bc \label{FINITE}
A finite para-tridiagonal (Hessenberg) matrix of order $N$ is unitary iff it has the
form $C(\bfa_N,\bfb_{N-1})$ ($H(\bfa_N,\bfb_{N-1})$), where
$|a_n|^2+|b_n|^2=1$ for $1 \leq n\leq N-1$ and $|a_N|=1$.
\ec

\bpr
This result is just a direct consequence of Theorem \ref{INFINITE} and the
following facts: a finite square matrix $M$ is unitary iff the infinite  matrix $M \oplus
I$ is unitary; the matrices $C(\bfa,\bfb)$ and $H(\bfa,\bfb)$, with $|a_n|^2+|b_n|^2=1,
\; \forall n\in\N$, decompose as a direct sum of their principal submatrices of order $N$
and an infinite matrix iff $b_N=0$.
\epr

\br \label{DECOMP}
{\it Decomposition property.}
Unitary para-tridiagonal and isometric Hessenberg matrices have similar decomposition
properties. They decompose as a sum of diagonal blocks iff, for some $N$, $b_N=0$, that
is, $a_N\in\T$. Moreover, in this situation, the blocks must again be unitary
para-tridiagonal and isometric Hessenberg matrices, respectively. Therefore, just looking
at the main diagonal we discover that, if $a_N\in\T$,
$$
\ba{l}
C(\bfa) = C(\bfa_N) \oplus C(\overline a_N\bfa^{(N)}),
\quad \bfa^{(N)}=(a_{N+n})_{n\in\N},
\\
C(\bfa_{N+M}) = C(\bfa_N) \oplus C(\overline a_N\bfa^{(N)}_M),
\quad \bfa^{(N)}_M=(a_{N+1},a_{N+2},\dots,a_{N+M}).
\ea
$$
Similar results hold for isometric Hessenberg matrices.
\er

\br \label{THETA}
{\it Factorization property.}
For $a,b\in\C$, let us define
$$
\Theta(a,b) := \pmatrix{-a & \overline b \cr b & \overline a},
\quad \hat\Theta_n(a,b) := I_{n-1} \oplus \Theta(a,b) \oplus I.
$$
Then, for any bounded sequences $\bfa,\bfb\in\C^{\aleph_0}$,
$$
C(\bfa,\bfb)=C_o(\bfa,\bfb)C_e(\bfa,\bfb)^T,
\quad H(\bfa,\bfb) = \prod_{n=1}^\infty \hat\Theta_n(a_n,b_n),
$$
where the infinite product, which has to be understood in the strong sense, is from
the left to the right, and
$$
C_e(\bfa,\bfb) = I_1 \oplus \left(\bigoplus_{n\in\N}
\Theta(a_{2n},b_{2n})\right),
\quad C_o(\bfa,\bfb) = \bigoplus_{n\in\N} \Theta(a_{2n-1},b_{2n-1}).
$$
These factorizations show explicitly the isometric properties of the matrices given
in Theorem \ref{INFINITE} and Corollary \ref{FINITE}.

We also denote
$\Theta(a):=\Theta(a,\sqrt{1-|a|^2})$,
$\hat\Theta_n(a) := I_{n-1} \oplus \Theta(a) \oplus I$, so that
$H(\bfa) = \prod_{n=1}^\infty \hat\Theta_n(a_n)$ and
$C(\bfa)=C_o(\bfa)C_e(\bfa)$,
where
$C_e(\bfa) := I_1 \oplus \left(\bigoplus_{n\in\N} \Theta(a_{2n})\right)$ and
$C_o(\bfa) := \bigoplus_{n\in\N} \Theta(a_{2n-1})$.
\er

In the case of finite unitary Hessenberg matrices, the above properties have been used
for spectral computations \cite{Gr86,GrRe90}. Notice that the factorization property in
the Hessenberg case is much worse than in the para-tridiagonal one.

We know that any unitary matrix represents an orthogonal sum of multiplication operators
on $\T$ and, hence, is unitarily equivalent to a direct sum of unitary irreducible
para-tridiagonal matrices. However, if the initial matrix is also para-tridiagonal,
the equivalence becomes an equality. This is just a consequence of Theorem
\ref{INFINITE}, Corollary \ref{FINITE} and the decomposition property given in Remark
\ref{DECOMP}. For the same reason, a similar result is also true for isometric
Hessenberg matrices.

\bc \label{RED}
Every unitary para-tridiagonal (isometric Hessenberg) matrix is a direct sum of
irreducible unitary para-tridiagonal (isometric Hessenberg) matrices.
\ec

Even more, in the study of irreducible unitary para-tridiagonal (isometric Hessenberg)
matrices, it is enough to consider those with the form $C(\bfa)$ ($H(\bfa)$) and their
principal submatrices. More precisely, we have the following immediate result.

\bl \label{AB-A}
For any $\bfa\in\C^{\aleph_0},\bfb\in(\C\backslash\{0\})^{\aleph_0}$ it is
$H(\bfa)=R^*H(\bfa,\bfb)R$ and $C(\bfa)=S^*C(\bfa,\bfb)S$, where
$$
\ba{l}
R=\pmatrix{r_1 & & & \cr & \hskip-3pt r_2 & & \cr
& & \hskip-3pt r_3 & \cr & & & \hskip-3pt \ddots},
\quad \ba{l} r_1=1, \\ r_{n+1}/r_n = b_n/|b_n|, \quad n\geq1, \ea
\\
S=\pmatrix{\overline s_1 & & & \cr & \hskip-3pt s_2 & & \cr
& & \hskip-3pt \overline s_3 & \cr & & & \hskip-3pt \ddots},
\quad \ba{l} s_1=1, \\ s_2=b_1/|b_1|, \\
s_{n+1}/s_{n-1} = \overline b_{n-1} b_n /|\overline b_{n-1} b_n|, \quad n\geq2. \ea
\ea
$$
If $R_N,S_N$ are the principal submatrices of order $N$ of $R,S$ respectively,
$H(\bfa_N)=R_N^*H(\bfa_N,\bfb_{N-1})R_N$ and
$C(\bfa_N)=S_N^*H(\bfa_N,\bfb_{N-1})S_N$.
\el

Notice that Theorem \ref{INFINITE} and the above lemma imply that an infinite
Hessenberg matrix is unitary iff it has the form $H(\bfa,\bfb)$ with
$|a_n|^2+|b_n|^2=1$ and $\bfa\notin\ell^2$.

The preceding results have the following consequence, that represents the
``five-diagonal reduction'' of the spectral problem for any unitary Hessenberg
matrix. Without loss of generality we consider only the irreducible case.

\bt \label{H-C}
Let $H=(h_{i,j})$ be a (finite or infinite) unitary irreducible Hessenberg matrix and
let us define
$$
\tau_n := \cases{
1, & $n=1$, \cr
\ds \prod_{k=1}^{n-1} h_{k+1,k}, & $n\geq2$.
}
$$
Then, $H$ is unitarily equivalent to a para-tridiagonal matrix $C=(c_{i,j})$ with the
form $C(\bfa)$ or $C(\bfa_N)$, where
$$
a_n = - {h_{1,n} \over \overline \tau_n}, \quad n\geq1.
$$
The unitary equivalence is given by $H=V^*CV$, where the columns of $V=(v_{i,j})$
can be recursively obtained by
$$
\ba{l}
v_{i,j} = \cases{
\ds {\overline \tau_j \over |\tau_j|} \, \delta_{i,j}, & $j=1,2$,
\cr
\ds {1 \over h_{j,j-1}}
\left( \sum_{k=i-2}^{\min\{i+2,2j-4\}} \kern-15pt c_{i,k} \, v_{k,j-1}
- \kern-15pt
\sum_{k=\min\{i,[{i+3\over2}]\}}^{j-1} \kern-15pt h_{k,j-1} \, v_{i,k} \right) \kern-3pt,
& $i \leq 2j-2$, $j \geq 3$,
\cr
0, & $i \geq 2j-1$, $j \geq 3$.
}
\ea
$$
In the above expression the sums have to be understood only over those terms in which
the matrix coefficients have indices between $1$ and the order of $H$.
The eigenvectors $x(\lambda) = \sum_n x_n(\lambda) e_n$ of $H$ and
$y(\lambda) = \sum_n y_n(\lambda) e_n$ of
$C$ corresponding to the same eigenvalue $\lambda$ are related by
$$
x_{2k-1}(\lambda) =
{\tau_{2k-1} \over |\tau_{2k-1}|} \lambda^{1-k} \overline{y_{2k-1}(\lambda)},
\quad x_{2k}(\lambda) =
{\tau_{2k} \over |\tau_{2k}|} \lambda^{1-k} y_{2k}(\lambda),
\quad k\geq1.
$$
\et

\bpr
We will consider only the case of an infinite matrix $H$, the proof for the finite
case being completely analogous. Then, from Theorem \ref{INFINITE}, $H$ must have the
form $H(\bfa,\bfb)$, $\bfa\in\D^{\aleph_0}\backslash\ell^2$. So, according to Lemma
\ref{AB-A}, $H$ is unitarily equivalent to $H(\bfa)$ which, at the same time, is
unitarily equivalent to $C(\bfa)$ since they are representations of the same
multiplication operator.

We know that $H=RH(\bfa)R^*$, $R=(r_i\delta_{i,j})$, where $r_i=\tau_i/|\tau_i|$ since
$h_{i+1,i}=b_i$. Besides, if $\mu$ is the measure related to the sequence $\bfa$ of Schur
parameters and $(\varphi_n)_{n\geq0}$, $(\chi_n)_{n\geq0}$ are the corresponding OP, OLP
respectively, then $H(\bfa) = U^* C(\bfa) U$, $U = (u_{i,j})$,
$u_{i,j} = \big< \varphi_{j-1},\chi_{i-1} \big>_\mu$.
Therefore, $H = V^* C(\bfa) V$, $V = U R^*$. For $j=1,2$, $\varphi_{j-1}=\chi_{j-1}$ and,
so, $v_{i,j} = \overline r_j\delta_{i,j}$. For the rest of the columns in $V$, if $j\geq2$,
$\varphi_j \in \spn\{1,z,z^{-1},\dots,z^{1-j},z^j\} = \spn\{\chi_0,\chi_1,\dots,\chi_{2j-1}\}$
and, thus, $\big< \varphi_j,\chi_i \big>_\mu = 0$ for $i \geq 2j$. Since $R$ is diagonal,
this implies that $v_{i,j} = 0$ for $i \geq 2j-1$, $j \geq 3$. Moreover, from the equality
$VH=C(\bfa)V$ we get
$$
\sum_{k=1}^j  v_{i,k} \, h_{k,j-1} =
\sum_{k=i-2}^{i+2} c_{i,k} \, v_{k,j-1}, \quad i \leq 2j-2, \; j \geq 3,
$$
which completes the expression given for $v_{i,j}$, once the restriction
$v_{i,j} = 0$, $i \geq \max \{j+1,2j-1\}$, is taken into account in the above sums.

Given an eigenvalue $\lambda$ of $H=RH(\bfa)R^*$, the corresponding eigenvectors are
spanned by
$\sum_{n\in\N} r_n \overline{\varphi_{n-1}(\lambda)}e_n$,
while the eigenvectors of $C(\bfa)$ are spanned by
$\sum_{n\in\N} \overline{\chi_{n-1}(\lambda)}e_n$.
Hence, the referred relation between eigenvectors is just a consequence of the
relation (\ref{OP-OLP}) between OP and OLP. \epr

The para-tridiagonal representations improve the Hessenberg representations of
unitary operators because of their greater simplicity. Besides, as we pointed out in
Remark \ref{DECOMP}, they have similar decomposition properties and, thus, ``Divide
and Conquer" algorithms \cite{GrRe90} can be also developed for the spectral problem
of a unitary para-tridiagonal matrix. Even more, the factorization given in Remark
\ref{THETA} allows to write the corresponding five-diagonal eigenvalue problem
equivalently as a generalized eigenvalue problem for a tri-diagonal pair of unitary
matrices.

Now we reach the announced result about the minimal representations of unitary
operators.

\bt \label{BAND}
A $(p,q)$-diagonal unitary matrix is a sum of diagonal blocks of order not greater than
$p+q$ if $p$ or $q$ are equal to 1.
\et

\bpr
We can restrict our attention to the case of $(1,q)$-diagonal matrices since, otherwise,
we can deal with the adjoint matrix, that keeps the unitarity. Also, it is enough to
prove that, if such a unitary band matrix has order greater than $q+1$, then it must
decompose as a sum of smaller diagonal blocks. Let us suppose that $\Omega$ is of order
greater that $q+1$ and does not decompose. $\Omega$ is, in particular, an isometric
Hesssenberg matrix and, hence, $\Omega=H(\bfa,\bfb)$ or $\Omega=H(\bfa_N,\bfb_{N-1})$,
$N \geq q+2$. Since it does not decompose, $b_n\neq0$ for all $n$. Thus, for $j \geq q+2$,
$\omega_{1,j} = - \prod_{k=1}^{j-1}\overline b_k a_j = 0$ implies $a_j=0$.
Hence, if  $i \leq j$,
$\omega_{i,j} = - \overline a_{i-1} a_j \prod_{k=i}^{j-1}\overline b_k = 0$ for
$j \geq q+2$. Therefore,
$\{\Omega^*e_1,\Omega^*e_2,\dots,\Omega^*e_{q+2}\} \subset
\spn\{e_1,e_2,\dots,e_{q+1}\}$,
which is a contradiction with the unitarity of $\Omega$.
\epr

A matrix that decomposes as a sum of finite diagonal blocks has always pure point
spectrum. Therefore, the previous theorem shows that the only $(p,q)$-diagonal
representations possible for any unitary operator are those where $p,q\ge2$.
Consequently, we have the following corollary.

\bc \label{minREP}
The minimal representations of unitary operators are five-diagonal.
\ec

\section{Krein's Theorem}

One of the advantages of band representations is that they make it easier to decide the
``smallness" of a perturbation.  For example, the compactness of an operator is
equivalent to stating that the diagonals of a band representation converge to 0. This
makes it quite simple, for example, to apply Weyl's Theorem \cite{We10,Ka66,ReSi78} for
the invariance of the essential spectrum. Also, it is easier to prove that a perturbation
belongs to the trace class, which can be used to give a simple and elegant operator
theoretic proof of Rakhmanov's lemma \cite{Si04} using the Kato-Rosenblum Theorem
\cite{Ka57,Ro57,Ka66,ReSi78} on the invariance of the absolutely continuous part of the
spectrum of an operator.

In spite of the difficulties that appear, many results about the orthogonality measures
of OP on $\T$ have been obtained using the Hessenberg representation
\cite{Go00a,Go00b,Go00c,GoNeAs95,Te92}, mainly due to the efforts of L. Golinskii. The
proofs of such results can be now simplified, but we want to show some new results and
advantages provided by the para-tridiagonal representation in the analysis of the
relation between measures and Schur parameters.

First of all, we will discuss the advantages found in the application of Krein's Theorem
\cite{AkKr38}, getting new results for discrete measures whose support has a finite
derived set. Krein's Theorem asserts that, given a measure $\mu$ with infinite support,
it is equivalent to saying that $\supp\mu$ accumulates on the finite set
$\{w_1,w_2,...,w_N\}$ and that the operator $p_N(U_\mu)$ is compact, where
$p_N(z)=\prod_{i=1}^N(z-w_i)$. This theorem was established in 1962 by N.I. Akhiezer and
M.G. Krein \cite{AkKr38} for measures  on the real line with finite moments. Recently,
the translation to the unit circle was given by L. Golinskii \cite{Go00a}, who succeeded
in characterizing in terms of the Schur parameters the measures whose support has one or
two limit points, using the Hessenberg representation of $U_\mu$. He also proved that the
Schur parameters of any measure on $\T$ whose support has a finite number $N$ of limit
points must satisfy
\beq \label{DISCRETE}
\lim_n \rho_n\rho_{n+1}\cdots\rho_{n+N-1} = 0,
\eeq
from which comes the property $\limsn|a_n|=1$ for any measure whose support has a
finite derived set.

However, with the Hessenberg representation it is hard to go further in this direction.
The para-tridiagonal representation makes things easier, not only because of its band
structure, but also due to its factorization properties. In the context of the
para-tridiagonal representation, for the application of Krein's Theorem it is necessary
to decide the compactness of $p_N(C(\bfa))$, where $\bfa$ is the sequence of Schur
parameters associated with $\mu$. This requires the calculation of the $4N+1$ diagonals
of $p_N(C(\bfa))$, some of them possibly giving redundant information. We can optimize
the calculations using the factorization of Remark \ref{THETA}.

\bp \label{PREVIO-KREIN}
Given $w_1,w_2,\dots,w_N\in\T$, let us define
$$
q_N(C(\bfa)) = \cases{
{C(\bfa)^*}^k p_N(C(\bfa)) & if $N=2k$,
\smallskip\cr
C_o(\bfa)^* {C(\bfa)^*}^k p_N(C(\bfa)) & if $N=2k+1$,
}
$$
where $p_N(z)=\prod_{i=1}^N(z-w_i)$.
Then, $q_N(C(\bfa))$ is a $2N+1$-diagonal matrix such that
$$
q_N(C(\bfa))^* = \cases{
\left(\prod_{i=1}^N\overline w_i\right) q_N(C(\bfa)), & if $N$ is even,
\smallskip\cr
\overline{q_N(C(\bfa))}, & if $N$ is odd,
}
$$
and $p_N(C(\bfa))$ is compact iff
$\lim_nq_N(C(\bfa))_{n+m,n}=0$ for $m=0,1,\dots,N$.
\ep

\bpr

From the unitarity of $C(\bfa)$ and $C_o(\bfa)$, the equivalence between the compactness
of $p_N(C(\bfa))$ and $q_N(C(\bfa))$ follows. The matrix $q_N(C(\bfa))$ is a linear
combination of products of, at most, $N$ tri-diagonal matrices, so, it is
$2N+1$-diagonal. Thus, $q_N(C(\bfa))$ is compact iff $\lim_nq_N(C(\bfa))_{n+m,n}=0$ for
$m=0,\pm1,\dots,\pm N$. Hence,  to finish the proof we just have to check the relations
between $q_N(C(\bfa))$ and $q_N(C(\bfa))^*$.

When $N$ is odd, $q_N(C(\bfa))$ is a linear combination of products of an odd number
of alternate factors $C_o(\bfa)$ and $C_e(\bfa)$, or their adjoints. Since $C_o(\bfa)$
and $C_e(\bfa)$ are symmetric, $q_N(C(\bfa))$ is symmetric too.

In the case of even $N=2k$, we can write
$q_N(C(\bfa))=\prod_{i=1}^k r_i(C(\bfa))$,
$r_i(C(\bfa))=(C(\bfa)-w_i)C(\bfa)^*(C(\bfa)-w_{k+i})$.
The result is just a consequence of the fact that
$r_i(C(\bfa))^*=\overline w_i\overline w_{k+i}r_i(C(\bfa))$.
\epr

Therefore, we can apply Krein's Theorem imposing only that the main and lower diagonals
of $q_N(C(\bfa))$ converge to 0, which will give in general $N+1$ asymptotic conditions
for the Schur parameters of a measure whose support has $N$ given limit points. For
illustrative purposes we present the results achieved using this procedure when applied
to the characterization of measures whose support has up to three limit points.

\bp \label{EX-KREIN}
Let $\mu$ be the measure associated with the sequence $\bfa$ of Schur parameters. Then:
\be
\item
$\{\supp\mu\}^{'}=\{\alpha\}$ iff

$\lim_n ( \overline a_na_{n+1}+\alpha )=0$.

\item
$\{\supp\mu\}^{'}\subset\{\alpha,\beta\}$ iff

$\lim_n \rho_n\rho_{n+1}=0$,

$\lim_n \rho_{n+1} ( \overline a_na_{n+2}-\alpha\beta )=0$,

$\lim_n ( \overline a_na_{n+1}
        +\alpha\beta a_n\overline a_{n+1}+\alpha+\beta )=0$.

\item
$\{\supp\mu\}^{'}\subset\{\alpha,\beta,\gamma\}$ iff

$\lim_n \rho_n\rho_{n+1}\rho_{n+2}=0$,

$\lim_n \rho_{n+1}\rho_{n+2}
        ( \overline a_na_{n+3}+\alpha\beta\gamma )=0$,

$\lim_n \rho_{n+1} ( \overline a_na_{n+1}+\overline a_{n+1}a_{n+2}
        -\alpha\beta\gamma a_n\overline a_{n+2}+\alpha+\beta+\gamma )=0$,

$\lim_n ( \overline a_na_{n+1}^2-\rho_{n+1}^2a_{n+2}
       +\alpha\beta\gamma(a_n^2\overline a_{n+1}-a_{n-1}\rho_n^2)\,+$
\vskip-3pt
$ \kern123.5pt +\,( \alpha\beta+\beta\gamma+\gamma\alpha)a_n
       +(\alpha+\beta+\gamma)a_{n+1} )=0$.
\ee
\ep

The first result of the above proposition is the same one obtained in \cite{Go00a},
but the second assertion simplifies the one given in \cite{Go00a}. Notice that the
relations given in the two last cases of the proposition also imply
$$
\lim_n \rho_{n+1} \left(\overline a_na_{n+1}+\overline a_{n+1}a_{n+2}
       +\alpha+\beta \right)=0,
$$
if $\{\supp\mu\}^{'}\subset\{\alpha,\beta\}$, while, for
$\{\supp\mu\}^{'}\subset\{\alpha,\beta,\gamma\}$,
$$
\lim_n \rho_{n+1}\rho_{n+2} \left( \overline a_na_{n+1}+\overline a_{n+1}a_{n+2}
       +\overline a_{n+2}a_{n+3}+\alpha+\beta+\gamma \right)=0.
$$

The above results suggest the following improvement of the property (\ref{DISCRETE}) that
gives a common feature for measures whose support has a finite derived set.

\bt \label{KREIN}
Let $\mu$ be the measure associated with the sequence $\bfa$ of Schur parameters.
If $\{\supp\mu\}^{'}\subset\{w_1,w_2,\dots,w_N\}$, then
$$
\ba{c}
\ds\lim_n \prod_{i=1}^N\rho_{n+i}=0, \quad N\ge1,
\\
\ds\lim_n \bigg(\overline a_n a_{n+N}-P\bigg)
          \prod_{i=1}^{N-1}\rho_{n+i}=0, \quad N\ge2,
\\
\ds\lim_n \left(\sum_{j=1}^N\overline a_{n+j-1}a_{n+j} + S\right)
          \prod_{i=1}^{N-1}\rho_{n+i}=0, \quad N\ge2,
\\
\ds\lim_n
\left(\sum_{j=1}^{N-1}\overline a_{n+j-1}a_{n+j} + Pa_n\overline a_{n+N-1} + S\right)
\prod_{i=1}^{N-2}\rho_{n+i}=0, \quad N\ge2,
\ea
$$
where $P:=(-1)^Nw_1w_2\cdots w_N$ and $S:=w_1+w_2+\cdots+w_N$.
\et

\bpr We will consider only the case of even $N=2k$, since the analysis for odd $N$ is
analogous. Then, the operator $q_N(C(\bfa))$ associated with the points
$w_1,w_2,\dots,w_N,$ given in Proposition \ref{PREVIO-KREIN}, has the form
$$
q_N(C(\bfa)) =
C(\bfa)^k-SC(\bfa)^{k-1}+\cdots-\overline SPC(\bfa)^*{^{k-1}}+PC(\bfa)^*{^k}.
$$
For $N-3 \leq m \leq N$, $q_N(C(\bfa))_{n+m,n}=q_N^{(1)}(C(\bfa))_{n+m,n}$, where
$$
q_N^{(1)}(C(\bfa)) =
C(\bfa)^k-SC(\bfa)^{k-1}-\overline SPC(\bfa)^*{^{k-1}}+PC(\bfa)^*{^k}.
$$
So, $\lim_nq_N^{(1)}(C(\bfa))_{n+m,n}=0$, $N-3 \leq m \leq N$, under the hypothesis
for $\mu$.

Let us examine the coefficients $q_N^{(1)}(C(\bfa))_{n+m,n}$ for $m=N,N-1,N-2$. We can
write $C_o(\bfa)=A_o+VB_o+B_oV^*$ and $C_e(\bfa)=A_e+VB_e+B_eV^*$, $V$ being the right
shift, defined by $Ve_n=e_{n+1}$, and
$$
\ba{l}
A_oe_n=\cases{-a_ne_n, & odd $n$, \cr \overline a_{n-1}e_n, & even $n$,} \quad
A_ee_n=\cases{\overline a_{n-1}e_n, & odd $n$, \cr -a_ne_n, & even $n$,}
\medskip\\
B_oe_n=\cases{\rho_ne_n, & odd $n$, \cr 0, & even $n$,} \quad \kern11pt
B_ee_n=\cases{0, & odd $n$, \cr \rho_ne_n, & even $n$.}
\ea
$$

Taking into account that $V$, $V^*$ rise and lower the indices of the vectors $e_n$,
respectively, $V^*e_1=0$, and $B_o$, $B_e$ vanish over vectors with even and odd index
$n$, respectively, we find that
$$
\ba{l}
(q_N^{(1)}(C(\bfa))e_n,e_{n+N}) =
(((VB_oVB_e)^k+P(VB_eVB_o)^k)e_n,e_{n+N}) =
\medskip\\
\kern105pt = \cases{
\rho_n\rho_{n+1}\cdots\rho_{n+N-1}, & even $n$,
\cr
P\rho_n\rho_{n+1}\cdots\rho_{n+N-1}, & odd $n$,
}
\ea
$$
$$
\ba{l}
(q_N^{(1)}(C(\bfa))e_n,e_{n+N-1}) =
\medskip\\
\kern30pt
= (((VB_oVB_e)^{k-1}VB_oA_e+A_oVB_e(VB_oVB_e)^{k-1})e_n,e_{n+N-1}) +
\medskip\\
\kern40pt + P (((VB_eVB_o)^{k-1}VB_eA_o^*+A_e^*VB_o(VB_eVB_o)^{k-1})e_n,e_{n+N-1}) =
\medskip\\
\kern115pt = \cases{
-\rho_n\cdots\rho_{n+N-2}(a_{n+N-1}-Pa_{n-1}), & even $n$,
\cr
-P\rho_n\cdots\rho_{n+N-2}(\overline a_{n+N-1}-\overline P\overline a_{n-1}), & odd $n$.
}
\ea
$$
Therefore, $\lim_n \prod_{i=1}^N\rho_{n+i}=0$ and $\lim_n
(a_{n+N}-Pa_n)\prod_{i=1}^{N-1}\rho_{n+i}=0$, which is equivalent to the first and second
equalities of the theorem.

Concerning the coefficients $(q_N^{(1)}(C(\bfa))e_n,e_{n+N-2})$, we have that
$$
\ba{l}
(C(\bfa)^ke_n,e_{n+N-2}) = ((VB_oVB_e\cdots VB_oVB_eA_oA_e)e_n,e_{n+N-2}) \, +
\medskip\\
\kern97pt + \, ((VB_oVB_e\cdots VB_oA_eA_oVB_e)e_n,e_{n+N-2}) + \cdots
\medskip\\
\kern80pt \cdots + ((A_oVB_e\cdots VB_oVB_eVB_oA_e)e_n,e_{n+N-2}) =
\medskip\\
\kern70pt = \cases{
-\rho_n\cdots\rho_{n+N-3}
(\overline a_{n-1}a_n+\cdots+\overline a_{n+N-3}a_{n+N-2}), & even $n$,
\cr
-\rho_n\cdots\rho_{n+N-3}\overline a_{n-1}a_{n+N-2}, & odd $n$,
}
\ea
$$
$$
\ba{l}
(C(\bfa)^*{^k}e_n,e_{n+N-2}) = ((VB_eVB_o\cdots VB_eVB_oA_e^*A_o^*)e_n,e_{n+N-2}) \, +
\medskip\\
\kern101pt + \,((VB_eVB_o\cdots VB_eA_o^*A_e^*VB_o)e_n,e_{n+N-2}) + \cdots
\medskip\\
\kern84pt \cdots + ((A_e^*VB_o\cdots VB_eVB_oVB_eA_o^*)e_n,e_{n+N-2}) =
\medskip\\
\kern70pt = \cases{
-\rho_n\cdots\rho_{n+N-3}a_{n-1}\overline a_{n+N-2}, & even $n$,
\cr
-\rho_n\cdots\rho_{n+N-3}
(a_{n-1}\overline a_n+\cdots+a_{n+N-3}\overline a_{n+N-2}), & odd $n$,
}
\ea
$$
$$
\ba{l} \kern-106pt
(C(\bfa)^{k-1}e_n,e_{n+N-2}) = ((VB_oVB_e)^{k-1}e_n,e_{n+N-2}) =
\medskip\\
= \cases{
\rho_n\rho_{n+1}\cdots\rho_{n+N-3}, & even $n$,
\cr
0, & odd $n$,
}
\ea
$$
$$
\ba{l} \kern-100pt
(C(\bfa)^*{^{k-1}}e_n,e_{n+N-2}) = ((VB_eVB_o)^{k-1}e_n,e_{n+N-2}) =
\medskip\\
\kern11pt = \cases{
0, & even $n$,
\cr
\rho_n\rho_{n+1}\cdots\rho_{n+N-3}, & odd $n$.
}
\ea
$$
From these results the last equality of the theorem follows. The third relation is just
a consequence of the other ones.
\epr

\section{Perturbations of the Schur parameters and mass points}

In the previous discussion we have exploited the band structure and factorization
properties of the para-tridiagonal representation. Now we will also show the advantages
of its simple dependence on the Schur parameters, in particular, of the fact that,
contrary to the Hessenberg representation, any Schur parameter appears in only a
finite number of elements of the para-tridiagonal representation.

The application to the study of OP of standard results of operator theory, like the Weyl,
Krein or Kato-Rosenblum theorems, gives information about the limit points of the support
of the orthogonality measure. However, for the analysis of isolated mass points other
tools are more appropriate. This last section illustrates the usefulness of the
para-tridiagonal representation for this purpose too. Our aim is to study the behaviour
of the isolated mass points of the measure under monoparametric perturbations of the
Schur parameters using the Hellmann-Feynman Theorem \cite{Fe39,He37,If88,IsZh88}.

\medskip

Let us suppose a sequence $\bfa(t)\in\D^{\aleph_0}$ depending on $t \in I$, where
$I$ is an interval of $\R$. A measure $\mu^t$ corresponds to each value of $t$.
The related OP and OLP sequences will be denoted by $\bsF^t:=(\varphi_n^t)_{n\geq0}$
and $\bfX^t:=(\chi_n^t)_{n\geq0}$ respectively.

Besides, let $u(t)$ be a function of $t \in I$ with values on $\T$. For each $t$ we can
consider the finitely supported measure $\mu_N^t$ corresponding to the parameters
$(a_1(t),\dots,a_{N-1}(t),u(t))$, whose finite segments of OP and OLP are respectively
$\bsF_N^t:=(\varphi_n^t)_{n=0}^{N-1}$ and $\bfX_N^t:=(\varphi_n^t)_{n=0}^{N-1}$. The
importance of such discrete measures is that they weakly converge to $\mu^t$ and, thus,
they provide the, so called, Szeg\H o quadrature formulas \cite{JoNjTh89} for the measure
$\mu^t$.

We are interested in the evolution with $t$ of the isolated mass
points of $\mu^t$, that is, the isolated eigenvalues of
$\cC(t):=C(\bfa(t))$. We will also analyze the movement of the
mass points of the discrete approximations $\mu_N^t$, that is, the
eigenvalues of $\cC_N(t):=C(\hat\bfa_N(t))$,
$\hat\bfa_N(t):=(a_1(t),\dots,a_{N-1}(t),u(t))$.

Since the finite matrices $\cC_N(t)$ have $N$ different eigenvalues, in any interval
where $\hat\bfa_N(t)$ is differentiable with respect to $t$, its eigenvalues are
differentiable functions $\lambda(t)$ \cite{Ka66}. Moreover, the corresponding
eigenvectors $\bfX_N^t(\lambda(t))$ are also differentiable in $t$, since $\bfX_N(z)$ is
a differentiable function of $a_1,\dots,a_{N-1},z$.

Concerning the infinite matrix $\cC(t)$, a similar result holds,
but only locally. More precisely, let us suppose that $\cC(t)$ is
differentiable in norm with bounded derivative $\cC^{'}\!(t)$ and
$\|\cC^{'}\!(t)\|$ locally bounded. Then, if $\lambda_0$ is an
isolated eigenvalue of $\cC(t_0)$, there exists a neighbourhood of
$t_0$ where $\cC(t)$ has an isolated eigenvalue $\lambda(t)$ which
is differentiable and such that $\lambda(t_0)=\lambda_0$.
Moreover, a related eigenvector can be chosen as a strongly
differentiable function of $t$ in a neighbourhood of $t_0$
\cite{If91}.

This last discussion justifies the following lemma.

\bl \label{BAND-DIFF}
Let $\Omega(t)=(\omega_{i,j}(t))_{i,j\in\N}$ be a bounded band matrix depending on a
parameter $t \in I$, $I$ being an interval of $\R$.
Assume that the coefficients $\omega_{i,j}(t)$ are twice differentiable and
$\sup_{i,j\in\N}|\omega_{i,j}^{'}(t)|$, $\sup_{i,j\in\N}|\omega_{i,j}^{''}(t)|$ are
locally bounded on $I$. Then, $\Omega(t)$ is differentiable in norm with bounded
derivative $\Omega^{'}\!(t):=(\omega_{i,j}^{'}(t))_{i,j\in\N}$ and
$\|\Omega^{'}\!(t)\|$ is locally bounded on $I$.
\el

\bpr
We can write
$\Omega(t)=\Omega_0(t)+\sum_{k=1}^N(V^k\Omega_k(t)+\Omega_{-k}(t){V^*}^k)$,
where $\Omega_k(t)$, $|k| \leq N$, are diagonal matrices and $V$
is the right shift. Hence, if the statement is true for the
matrices $\Omega_k(t)$, it is also true for $\Omega(t)$. So, we
just have to check the proposition for a diagonal matrix
$\Omega(t)=\mbox{diag}(\omega_1(t),\omega_2(t),\dots)$. If
$\omega_n(t)$ are differentiable and
$\sup_{n\in\N}|\omega_n^{'}(t)|$ is locally bounded on $I$,
$\Omega^{'}\!(t):=\mbox{diag}(\omega_1^{'}(t),\omega_2^{'}(t),\dots)$
is bounded with
$\|\Omega^{'}\!(t)\|=\sup_{n\in\N}|\omega_n^{'}(t)|$ locally
bounded on $I$. If, besides, $\omega_n(t)$ are twice
differentiable and $\sup_{n\in\N}|\omega_n^{''}(t)| \leq K$ in a
neighbourhood of $t_0$, using the mean value theorem we get
$$
\left\|{\Omega(t)-\Omega(t_0) \over t-t_0}-\Omega^{'}(t_0)\right\| =
\sup_{n\in\N} \left|{\omega_n(t)-\omega_n(t_0) \over t-t_0}-\omega_n^{'}(t_0)\right|
\leq K|t-t_0|,
$$
for $t$ in such a neighbourhood. This proves the differentiability in norm.
\epr

Now we can state the following result for a differentiable monoparametric perturbation
of the Schur parameters.

\bp \label{HF1}
Let $a_n \colon I\to\D$ be differentiable for $n\in\N$. Then:
\be
\item
If $u \colon I\to\T$ is differentiable, the mass points of $\mu_N^t$
are differentiable functions $\lambda \colon I\to\T$ satisfying
$$
\lambda^{'}\!(t) = \mu_N^t(\{\lambda(t)\}) \,
X_N^t(\lambda(t))^T \cC_N^{'}\kern-0.5pt(t) \,\overline{X_N^t(\lambda(t))}.
$$
\item
If $a_n \colon I \to \D$ is twice differentiable for $n\in\N$ and
$\sup_{n\in\N}|a_n^{'}(t)|$, $\sup_{n\in\N}|a_n^{''}(t)|$,
$\sup_{n\in\N}|\rho_n^{'}(t)|$, $\sup_{n\in\N}|\rho_n^{''}(t)|$
are locally bounded on $I$, for any isolated mass point
$\lambda_0$ of $\mu^{t_0}$ there exists a differentiable function
$\lambda \colon J\to\T$ on a neighbourhood $J$ of $t_0$ such that
$\lambda(t)$ is an isolated mass point of $\mu^t$ for $t \in J$
and $\lambda(t_0)=\lambda_0$. This function satisfies
$$
\lambda^{'}\!(t) = \mu^t(\{\lambda(t)\}) \,
X^t(\lambda(t))^T \cC^{'}\!(t) \,\overline{X^t(\lambda(t))}.
$$
\ee
\ep

\bpr From the previous discussions and Lemma \ref{BAND-DIFF} we
find that the referred differentiable functions $\lambda(t)$ exist
under the conditions of the theorem. The expression for
$\lambda^{'}\!(t)$ follows from the Hellmann-Feynman Theorem for
normal operators. Let us consider the infinite case since the
analysis of the finite case is analogous. The mass points
$\lambda(t)$ are simple eigenvalues of $\cC(t)$ with associated
eigenspace spanned by $\overline{X^t(\lambda(t))}$. We know that
there exists a strongly differentiable eigenvector $Y(t)$ of
$\cC(t)$ with respect to $\lambda(t)$. Therefore, just
differentiating the equality
$Y(t)^*\cC(t)Y(t)=\lambda(t)Y(t)^*Y(t)$ and bearing in mind the
unitarity of $\cC(t)$, we get
$Y(t)^*\cC^{'}\!(t)Y(t)=\lambda^{'}\!(t)Y(t)^*Y(t)$. This relation
is also true when substituting $Y(t)$ by
$\overline{X^t(\lambda(t))}$ since they are proportional. The
statement 2 is then a consequence of the equality
$X^t(\lambda(t))^T \overline{X^t(\lambda(t))} =
\sum_{n\in\N}|\varphi_n(\lambda(t))|^2 = 1/\mu^t(\{\lambda(t)\})$.
\epr

It is natural to expect a qualitatively different behaviour of the
measure under rotations or dilatations of the Schur parameters.
Therefore, it could be interesting to examine the preceding result
when we decompose the monoparametric perturbation in the way
$a_n(t)=r_n(t)e^{i\alpha_n(t)}$, $r_n(t),\alpha_n(t)$ being real
functions.

\bt \label{HF2}
Let $a_n(t)=r_n(t)e^{i\alpha_n(t)}$, where $r_n \colon I\to(-1,1)$,
$\alpha_n \colon I\to\R$ are differentiable for $n\in\N$.
We define the functions
$$
\ba{l}
\ds\Gamma_n^t(z):={2\over \rho_n(t)^2} \,
\im\!\left(e^{-i\alpha_n(t)}z^{2-n}(\varphi_{n-1}^t(z))^2\right), \quad n\in\N,
\medskip\\
\Delta_n^t(z):=|\varphi_{n-1}^t(z)|^2-|\varphi_n^t(z)|^2, \quad n\in\N.
\ea
$$
\be
\item
If $u(t)=e^{i\beta(t)}$, being $\beta \colon I\to\R$ differentiable, the mass points of
$\mu_N^t$ have the form $\lambda(t)=e^{i\theta(t)}$, where $\theta \colon I\to\R$ is a
differentiable function that satisfies
$$
\ba{l}\,
\ds\theta^{'}\!(t) = \mu_N^t(\{\lambda(t)\}) \bigg\{ \sum_{n=1}^{N-1}
\left( r_n^{'}(t)\Gamma_n^t(\lambda(t))+\alpha_n^{'}(t)\Delta_n^t(\lambda(t))\right) +
\\ \kern240pt
+ \beta^{'}\!(t) |\varphi_{N-1}(\lambda(t))|^2 \bigg\}.
\ea
$$
\item
Let $r_n \colon I \to (-1,1)$, $\alpha_n \colon I \to \R$ be twice
differentiable. Assume that $\sup_{n\in\N}|r_n^{'}(t)|$,
$\sup_{n\in\N}|r_n^{''}(t)|$, $\sup_{n\in\N}|\alpha_n^{'}(t)|$,
$\sup_{n\in\N}|\alpha_n^{''}(t)|$ are locally bounded on $I$ and
there exists $r<1$ such that $|r_n(t)|<r$ whenever
$r_n^{'}(t)\neq0$. Then, if $\lambda_0$ is an isolated mass point
of $\mu^{t_0}$, there exists a differentiable function $\lambda
\colon J\to\T$ on a neighbourhood $J$ of $t_0$ such that
$\lambda(t)$ is an isolated mass point of $\mu^t$ for $t \in J$
and $\lambda(t_0)=\lambda_0$. $\lambda(t)=e^{i\theta(t)}$,
where $\theta \colon I\to\R$ is a differentiable function that
satisfies
$$
\ba{l}\,
\ds\theta^{'}\!(t) = \mu^t(\{\lambda(t)\}) \sum_{n=1}^\infty
\left( r_n^{'}(t)\Gamma_n^t(\lambda(t))+\alpha_n^{'}(t)\Delta_n^t(\lambda(t))\right).
\ea
$$
\ee
\et

\br \label{SERIE}
Notice that the above series converges due to the suppositions about the sequences
$(r_n(t))_{n\in\N}$, $(\alpha_n(t))_{n\in\N}$ and the fact that
$\bsF^t(\lambda(t))\in\ell^2$ since $\lambda(t)$ is a mass point of $\mu^t$.
\er

\bpr
The conditions given for the sequences $(r_n(t))_{n\in\N}$, $(\alpha_n(t))_{n\in\N}$
and the function $\beta(t)$ are enough to apply Proposition \ref{HF1}.
Therefore, the referred differentiable functions $\lambda(t)$ exist. Since $\R$ is
the universal covering space of $\T$, with the imaginary exponential as a covering map,
there exists a unique continuous real valued function $\theta(t)$ such that
$\lambda(t)=e^{i\theta(t)},\theta(t_0)=\mbox{Arg}(\lambda(t_0))$.
Moreover, the imaginary exponential is locally invertible with differentiable inverse,
so, $\theta(t)$ must be differentiable too. From Proposition \ref{HF1} we know that
$$
\theta^{'}\!(t) = {\lambda^{'}\!(t) \over i\lambda(t)} =
{\mu^t(\{\lambda(t)\}) \over i\lambda(t)} \,
X^t(\lambda(t))^T \cC^{'}\!(t) \,\overline{X^t(\lambda(t))}.
$$

The rest of the proof is just the calculation of the right hand side of the above
expression, which we will do only for the infinite case since the arguments in finite
case are similar. We can easily do this calculation using the factorization
$\cC(t)=\cC_o(t)\cC_e(t)$, $\cC_o(t)=C_o(\bfa(t))$, $\cC_e(t)=C_e(\bfa(t))$,
given by Remark \ref{THETA}. As a consequence of (\ref{RR-OP}) and (\ref{OP-OLP}),
$\cC_e(t) \overline{X^t(\lambda(t))} = X^t(\lambda(t))$ and
$\cC_o(t) X^t(\lambda(t)) = \lambda(t) \overline{X^t(\lambda(t))}$. Therefore,
$$
X^t(\lambda(t))^T \cC^{'}\!(t) \,\overline{X^t(\lambda(t))} =
X^t(\lambda(t))^T \cC_o^{'}(t) \, X^t(\lambda(t)) \, +
$$
\vskip-20pt
$$
\kern150pt + \lambda(t) X^t(\lambda(t))^* \cC_e^{'}(t) \,\overline{X^t(\lambda(t))},
$$
which, using (\ref{OP-OLP}), gives
$$
\ba{l}
\ds \theta^{'}\!(t) =
i\mu^t(\{\lambda(t)\}) \, \sum_{n=1}^\infty \lambda(t)^{-n}
\bigg\{ a_n^{'}(t) \varphi_{n-1}^{t*}(\lambda(t))^2
- \overline{a}_n^{'}(t) \varphi_n^t(\lambda(t))^2 -
\\  \kern230pt
- 2\rho_n^{'}(t) \varphi_{n-1}^{t*}(\lambda(t)) \varphi_n^t(\lambda(t)) \bigg\}.
\ea
$$
Finally, the expression given in the theorem for $\theta^{'}(t)$ follows from the above
one, taking into account (\ref{RR-OP}) and the relations
$$
a_n^{'}(t) = r_n^{'}(t) e^{i\alpha_n(t)} + i\alpha_n^{'}(t) a_n(t),
\quad \rho_n^{'}(t) = - {r_n(t) \over \rho_n(t)} \, r_n^{'}(t).
$$
\epr

From the above theorem we directly get a bound for the angular velocity of
the isolated mass points.

\bc \label{BOUND}
Under the conditions of Theorem \ref{HF2}
$$
|\theta^{'}(t)| \leq {2 \over 1-r^2} \sup_{n\geq1}|r_n^{'}(t)| +
\sup_{n\geq1}|\alpha_n^{'}(t)-\alpha_{n-1}^{'}(t)|,
$$
where $\alpha_0=0$ and, in the case of $\mu_N^t$, the sums run from 1 to $N$ and
$\alpha_N=\beta$.
\ec

The particular case of uniform rotations of the Schur parameters is specially
interesting. It has been previously considered in \cite{GoNe01} and by the authors in
\cite{Ca97,MPOP,FIVE}.

\bc \label{ROT}
Let $\bfa\in\D^{\aleph_0}$, $u\in\T$ and $\alpha \colon I\to\R$ differentiable.
If $a_n(t)=a_ne^{i\alpha(t)}$ for $n\in\N$ and $u(t)=ue^{i\alpha(t)}$, then:
\be
\item
The differentiable arguments of the mass points of $\mu_N^t$ satisfy
$$
\ba{l}\,
\ds\theta^{'}\!(t) = \mu_N^t(\{\lambda(t)\}) \, \alpha^{'}\!(t) .
\ea
$$
\item If $\alpha^{''}(t)$ exists and is locally bounded on $I$, the differentiable
arguments of the isolated mass points of $\mu^t$ satisfy
$$
\ba{l}\,
\ds\theta^{'}\!(t) = \mu^t(\{\lambda(t)\}) \, \alpha^{'}\!(t) .
\ea
$$
\ee
\ec

\bpr
Apply Theorem \ref{HF2} to $r_n(t)=|a_n|$, $\alpha_n(t)=\alpha(t)+\mbox{Arg}(a_n)$ and
$\beta(t)=\alpha(t)+\mbox{Arg}(u)$. Notice that $\alpha^{'}(t)$ is locally bounded on
$I$ if $\alpha(t)$ is twice differentiable.
\epr

This result is the generalization to arbitrary measures of the one founded in \cite{MPOP}
for finitely supported measures using the Hessenberg representation. It says that under a
uniform rotation of the Schur parameters the isolated mass points of the corresponding
measure rotate in the same direction and the  mass of each point gives its relative
angular velocity with respect to the angular velocity of the Schur parameters. Therefore,
a mass point rotates so much more quickly with the Schur parameters as its mass gets
bigger. In fact, Theorem \ref{HF2} suggests that, in general, the mass of an isolated
mass point gives a measure of its instability under perturbations of the Schur
parameters.

The study of the relation between Schur parameters and measures implies the attempt to
find families of Schur parameters associated with measures with some common features.
Theorem \ref{HF2} opens a way to find monoparametric families of Schur parameters whose
measures have a common mass point. Among the ways to do this, we will just select some of
them.

\subsection{Measures with a fixed mass point}

Let $\mu$ be the measure corresponding to a sequence $\bfa=(r_ne^{i\alpha_n})_{n\in\N}$
of Schur parameters and $(\varphi_n)_{n\ge0}$ the associated OP. If $\lambda=e^{i\theta}$
is an isolated mass point of $\mu$, our aim is to find monoparametric perturbations
$\bfa(t)$, $\bfa(t_0)=\bfa$, such that the corresponding measures $\mu^t$ have the same
mass point, at least in a neighbourhood of $t_0$. We will also consider the analogous
problem for the finitely supported measures $\mu_N$ associated with the parameters
$(a_1,\dots,a_{N-1},u)$, $u=e^{i\beta}$. In what follows we suppose that the perturbation
satisfies the conditions given in Theorem \ref{HF2}.

\bigskip
\noindent {\bf Case 1.}
$a_k(t)=\cases{a_k & if $k \neq n$, \cr r(t)e^{i\alpha(t)} & if $k=n$.}
\qquad (r(t_0)=r_n, \alpha(t_0)=\alpha_n)$
\medskip

This case corresponds to the perturbation of only the $n$-th Schur parameter. So, the
first $n$ OP coincide with the unperturbed ones. Using (\ref{RR-OP}) we get from Theorem
\ref{HF2} that $\lambda=e^{i\theta}$ is a fixed mass point if
$$
\ba{l}
r^{'}\!(t) \,
\im\!\left(e^{-i\alpha(t)}\lambda^{2-n}(\varphi_{n-1}(\lambda))^2\right) =
\medskip \\ \kern90pt
= \alpha^{'}\!(t) r(t)
\left\{ \re\!\left(e^{-i\alpha(t)}\lambda^{2-n}(\varphi_{n-1}(\lambda))\right) +
r(t)^2 |\varphi_{n-1}(\lambda)|^2 \right\}.
\ea
$$
If $(\varphi_{n-1}(\lambda))^2=|\varphi_k(\lambda)|^2e^{i\xi}$, the above equation becomes
$$
d(r \sin(\alpha+(n-2)\theta-\xi)) + r^2 d\alpha =0,
$$
whose solution for the conditions $r(t_0)=r_n,\alpha(t_0)=\alpha_n$ is
$$
\ba{l}
\ds r = {\sin c \over \sin(\alpha+(n-2)\theta-\xi-c)},
\medskip\\
\ds c = \arctan\left({r_n\sin(\alpha_n+(n-2)\theta-\xi) \over
1+r_n\cos(\alpha_n+(n-2)\theta-\xi)}\right).
\ea
$$
The same solution appears in the case of a finitely supported measure $\mu_N$, $N>n$,
if we leave the parameter $u$ unperturbed.

\bigskip
\noindent {\bf Case 2.}
$a_k(t)=\cases{a_k & if $k<n$, \cr r(t)e^{i\alpha(t)} & if $k=n$, \cr
e^{i(\alpha(t)-\alpha_n)}a_k & if $k>n$.}
\qquad (r(t_0)=r_n, \alpha(t_0)=\alpha_n)$
\medskip

Again, the first $n$ OP coincide with the unperturbed ones. The condition given by
Theorem \ref{HF2} for a fixed mass point $\lambda=e^{i\theta}$ is now
$$
\ba{l}
\ds {2dr \over 1-r^2} \sin(\alpha+(n-2)\theta-\xi) - d\alpha = 0,
\ea
$$
where $\xi$ is the phase of $(\varphi_{n-1}(\lambda))^2$.
The solution for the conditions $r(t_0)=r_n,\alpha(t_0)=\alpha_n$ is
$$
\ba{l}
\ds r = {\sin{1 \over 2}(\alpha+(n-2)\theta-\xi-c) \over
\sin{1 \over 2}(\alpha+(n-2)\theta-\xi+c)},
\medskip\\
\ds c = 2\arctan\left( {1-r_n \over 1+r_n}
\tan{1 \over 2}(\alpha_n+(n-2)\theta-\xi) \right).
\ea
$$
This solution remains valid in the case of a measure $\mu_N$, $N>n$, if we also include
a perturbation $u(t)=e^{i(\alpha(t)-\alpha_n)}u$ of the parameter $u$.

\bigskip
\noindent {\bf Case 3.}
$a_k(t)=\cases{a_k & if $k<n$, \cr r(t)e^{i\alpha_n} & if $k=n$, \cr
e^{i\alpha(t)}a_k & if $k>n$.}
\qquad (r(t_0)=r_n, \alpha(t_0)=0)$
\medskip

As in the previous cases, the first $n$ OP coincide with the unperturbed ones.
From Theorem \ref{HF2} and using (\ref{RR-OP}) we find that the perturbations of this
type with a fixed mass point $\lambda=e^{i\theta}$ are characterized by
$$
\ba{l}
\ds 2\sin(\alpha_n+(n-2)\theta-\xi) \, dr
- (1+2r\cos(\alpha_n+(n-2)\theta-\xi)+r^2) \, d\alpha = 0,
\ea
$$
where, again, $\xi$ is the phase of $(\varphi_{n-1}(\lambda))^2$.
The solution for the conditions $r(t_0)=r_n,\alpha(t_0)=0$ is
$$
\ba{l}
\ds r = -{\sin({1 \over 2}\alpha+\alpha_n+(n-2)\theta-\xi-c) \over
\sin({1 \over 2}\alpha-c)},
\medskip\\
\ds c = \arctan \left(
{\sin(\alpha_n+(n-2)\theta-\xi) \over
r_n+\cos(\alpha_n+(n-2)\theta-\xi)}\right),
\ea
$$
which is also valid in the case of a measure $\mu_N$, $N>n$, if including a perturbation
$u(t)=e^{i\alpha(t)}u$ of $u$.

\medskip

If $\mu^{t_0}$ has an isolated point at $\lambda=e^{i\theta}$, the previous relations
between $r$ and $\alpha$ provide perturbations of the Schur parameters that give families
of measures with the same mass point $\lambda$, at least for $r$, $\alpha$ in a
neighbourhood of $r(t_0)$, $\alpha(t_0)$. In the case of a finitely supported measure
this neighbourhood is only restricted by the condition $|r|<1$.

The simplest case of the above perturbations happens when $n=1$, where always
$\xi=0$. Another particularly simple situation is the perturbation of a Geronimus
measure, that corresponds to a constant sequence of Schur parameters.

\medskip
\noindent {\bf \underline{Example}:} perturbations of Geronimus measures with a fixed
mass point.
\smallskip

Let us consider the measure corresponding to a constant sequence of Schur parameters
$a_n=a\in\D\backslash\{0\}$, $n\ge1$ \cite{Ge61}. This measure has an isolated mass point
at $\lambda=(1-a)/(1-\overline a)$ if $|a-1/2|>1/2$, that is, if $ \re(a)<|a|^2$. The
related orthogonal polynomials are $\varphi_n(z)=\rho^{-n}(u_{n+1}(z)-(1-a)u_n(z))$,
where $\rho=\sqrt{1-|a|^2}$, $u_n(z)=(w_1(z)^n-w_2(z)^n)/(w_1(z)-w_2(z))$ and $w_1(z)$,
$w_2(z)$ are the solutions of $w^2-(z+1)w+\rho^2z=0$ \cite{Go99}. If
$\lambda=e^{i\theta}$, the phase of $w_1(\lambda)$ and $w_2(\lambda)$ is $\theta/2$ and,
thus, $(n-1)\theta/2$ is the phase of $u_n(\lambda)$. Hence, $\xi=(n-1)\theta$ and
$\alpha+(n-2)\theta-\xi=\alpha-\theta$ for all $n$. In this case, the relations between
$r$ and $\alpha$ that give a fixed mass point for the three previous perturbations are
independent of the index $n$ of the Schur parameter where the perturbation starts.

Let us write $a=r_0e^{i\alpha_0}$, $r_0,\alpha_0\in\R$. Then, the condition for the
existence of an isolated mass point is $\cos\alpha_0<r_0$. Using the explicit form of the
mass point we find that
$$
\cos(\alpha_0-\theta) = {(1+r_0^2)\cos\alpha_0-2r_0 \over 1+r_0^2-2r_0\cos\alpha_0},
\quad
\sin(\alpha_0-\theta) = {(1-r_0^2)\sin\alpha_0\over 1+r_0^2-2r_0\cos\alpha_0}.
$$
Taking into account these expressions we can find explicitly the relations between $r$
and $\alpha$ that give a fixed mass point at $\lambda=(1-a)/(1-\overline a)$ in the case
of the three perturbations previously studied. We find the following results:

\bigskip

\noindent Case 1.
$\ds r = {r_0 \sin\alpha_0 \over \sin\alpha-r_0\sin(\alpha-\alpha_0)}$.

\bigskip

\noindent Case 2.
$\ds r = {r_0 \sin{1\over2}(3\alpha_0-\alpha) + \sin{1\over2}(\alpha-\alpha_0)
\over \sin{1\over2}(\alpha_0+\alpha) - r_0\sin{1\over2}(\alpha-\alpha_0)}$.

\bigskip

\noindent Case 3.
$\ds r = {r_0 \sin{1\over2}(2\alpha_0-\alpha) + \sin{1\over2}\alpha
\over \sin{1\over2}(2\alpha_0-\alpha) + r_0\sin{1\over2}\alpha}$.

\bigskip

Notice that in the first and second cases $r=r_0$ for $\alpha=\alpha_0$, due to the
initial conditions $r(t_0)=r_0, \alpha(t_0)=\alpha_0$, while in the third case $r=r_0$
for $\alpha=0$, since $r(t_0)=r_0,\alpha(t_0)=0$. If $t_0=0$, we can choose
$\alpha(t)=\alpha_0+t$ in the first two cases and $\alpha(t)=t$ in the third one. Then,
as a consequence of the previous results, we find that, in a neighbourhood of $t=0$, the
following families $\bfa(t)$ of Schur parameters are related to measures with a common
mass point at $\lambda=(1-a)/(1-\overline a)$ (we assume $a\in\D\backslash\R$ and
$\re(a)<|a|^2$):

\bigskip

\noindent Case 1.
$a_k(t)=\cases{
a, & if $k \neq n$,
\cr
\ds{\im(a) \over \im(a)\cos t-(|a|^2-\re(a))\sin t}\,e^{it}a, & if $k=n$.}$

\bigskip

\noindent Case 2.
$a_k(t)=\cases{
a, & if $k<n$,
\cr
\ds{\im(a)\cos{t \over 2}+(1-\re(a))\sin{t \over 2} \over
\im(a)\cos{t \over 2}-(|a|^2-\re(a))\sin{t \over 2}}\,e^{it}a, & if $k=n$,
\cr
e^{it}a, & if $k>n$.}$

\bigskip

\noindent Case 3.
$a_k(t)=\cases{
a, & if $k<n$,
\cr
\ds{\im(a)\cos{t \over 2}+(1-\re(a))\sin{t \over 2} \over
\im(a)\cos{t \over 2}+(|a|^2-\re(a))\sin{t \over 2}}\,a, & if $k=n$,
\cr
e^{it}a, & if $k>n$.}$

\bigskip

\bigskip
\noindent{\bf Acknowledgements}
\medskip

The work of the authors was supported by Project E-12/25 of DGA (Diputaci\'on General de
Arag\'on) and by Ibercaja under grant IBE2002-CIEN-07.

\end{document}